\documentclass[a4paper,11pt]{article}
\usepackage{epsfig}
\usepackage{epstopdf}
\usepackage{amsmath, amssymb, graphics, todonotes, algorithm, algpseudocode, bm}

\catcode`@=11
\@addtoreset{equation}{section}
\catcode`@=12
\newcommand{\R}{{\mathbb R}}
\newcommand{\N}{\mathbb N}
\renewcommand{\P}{\mathbb P}
\newcommand{\Q}{\mathbb Q}

\newcommand{\I}{\mathbb I}

\DeclareMathOperator*{\argmin}{arg\,min}

\newtheorem{Theorem}{Theorem}

\begin{document}
\title{A semi-Lagrangian algorithm in policy space\\ for hybrid optimal control problems}
\author{Roberto Ferretti \footnote{Dipartimento di Matematica e Fisica,
Universit\`a Roma Tre, L.go S. Leonardo Murialdo, 1, 00146 Roma (Italy), e-mail: {\tt ferretti@mat.uniroma3.it}} 
\and Achille Sassi \footnote{Unit\'e de Math\'ematiques Appliqu\'ees -- ENSTA Paristech, 828 Boulevard des Mar\'echaux, 91120 Palaiseau (France), e-mail: {\tt ach.sassi@gmail.com}}}

\maketitle
\begin{abstract}
The mathematical framework of hybrid system is a recent and general tool to treat control systems involving control action of heterogeneous nature. In this paper, we construct and test a semi-Lagrangian numerical scheme for solving the Dynamic Programming equation of an infinite horizon optimal control problem for hybrid systems. In order to speed up convergence, we also propose an acceleration technique based on policy iteration. Finally, we validate the approach via some numerical tests in low dimension.
\end{abstract}
\vskip1cm\noindent

\section{Introduction}

In the last two decades, the concept of {\em hybrid control system} has provided a sound mathematical framework for treating control systems in which continuous and discrete control actions mix together, and this framework has also been successfully adapted to optimal control problems. Among the various systems covered by this theory, we mention economic models with restocking, multigear and hybrid vehicles, and, more in general, systems with switchings in the dynamics and/or impulsive changes in the state. In this work, we study efficient numerical methods for applying Dynamic Programming techniques to hybrid optimal control problems of infinite horizon type.

Among the various mathematical formulations of optimal control problems for hybrid systems, we will adopt here the one given in \cite{BraBorMit, DhaRam}. Let $\I=\{1,\ldots,m\}$, and consider the controlled system $(X,Q)$ satisfying:
\begin{equation} \label{eq_stato}
 \begin{cases}
 \dot{X}(t)=f\big(X(t),Q(t),\alpha(t)\big) \quad t \in (0,\infty)\\
 X(0)= x\\
 Q(0^+)=q
 \end{cases}
\end{equation}
where $x\in \R^d $, and $q\in \I$. Here, $X$ and $Q$ denote respectively the continuous and the discrete component of the state. The function $f :\R^d\times \I\times U \to \R^d$ is the continuous dynamics, for a set of continuous controls given by:
\[
\mathcal{U}=\big\{\alpha:(0,\infty) \to U \> | \> \alpha \text{ measurable}, \ U \mbox{ compact} \big\}.
\]
The trajectory undergoes discrete transitions when it enters two predefined sets $A$ (the {\em autonomous} jump set) and $C$ (the {\em controlled} jump set), both of them subsets of $\R^d\times\I$. More precisely:

\begin{itemize}

 \item On hitting $A$, the trajectory jumps to a predefined destination set $D\subset \R^d \times \I$. The jump driven by a transition map $g: \R^d \times \I \times \mathcal{V} \to D$, where $\mathcal{V}$ is a discrete finite control set. Denoting by $\tau_i$ a time at which the trajectory hits $A$, the new state will be $\big(X(\tau_i^+\big), Q(\tau_i^+)\big)=g\big(X(\tau_i^-),Q(\tau_i^-),w_i\big)$, for a control $w_i\in \mathcal{V}$.

 \item Entering the set $C$, the controller can choose either to jump or not. If the controller chooses to jump, then the continuous trajectory is moved to a new point in $D$. Denoting by $\xi_k$ one such time of jump, we will have $\big(X(\xi_k^-), Q(\xi_k^-)\big) \in  C$ and $(x',q')=\big(X(\xi_k^+), Q(\xi_k^+)\big) \in  D$.

\end{itemize}

The trajectory starting from $x\in \R^d$ with discrete state $q\in \I$ is therefore composed of branches of continuous evolution given by \eqref{eq_stato} between two discrete jumps at the transition times $\tau_i$ or $\xi_k$.

Now, considering an optimal control problem in the infinite horizon form, and including all control actions in a control strategy
$$
\theta:=\Big(\alpha(\cdot),\{w_i\}_{i\in\N},\big\{(\xi_k,x'_k,q'_k)\big\}_{k\in\N}\Big)
$$
we associate to $\theta$ a cost defined by:
\begin{equation}\label{J}
\begin{aligned}
 J(x,q;\theta) := & \int_0^{+\infty} \ell\big(X(t),Q(t),\alpha(t)\big)e^{-\lambda t} dt  \\
 & + \sum_{i=0}^\infty c_A\big(X(\tau_i^-),Q(\tau_i^-),w_i\big)e^{-\lambda \tau_i}  \\
 & + \sum_{k=0}^\infty c_C\big(X(\xi_k^-),Q(\xi_k^-),X(\xi_k^+),Q(\xi_k^+)\big)  e^{-\lambda \xi_k}
\end{aligned}
\end{equation}
where $\lambda>0$ is the discount factor, $\ell: \R^d\times \I \times U \to \R_+$ is the running cost, $c_A: A\times \mathcal{V} \to \R_+$ is the autonomous transition cost and $c_C:C\times D \to \R_+$ is the controlled transition cost. The value function $V$ of the problem is then defined as:
\begin{eqnarray} \label{f_valore}
 V(x,q):= \inf_\theta J(x,q;\theta).
\end{eqnarray}
We point out that, in this generality, the problem requires strong assumptions to be mathematically well-posed. In particular, it should be ensured that the value function \eqref{f_valore} is continuous, and that the so-called ``Zeno executions'' (i.e., the occurrence of an infinite number of transitions in a finite time interval) are avoided. We will give in the next section a precise set of assumptions, whereas in the examples we will apply the numerical technique under consideration in more general situations, showing that the recipe is robust enough to handle them.

To the best of our knowledge, the first rigorous theoretical study of the convergence of numerical schemes for the approximation of the value function of \eqref{eq_stato}--\eqref{J} has been given in \cite{FerZid}. Here, solvability of the scheme by value iteration is proved, along with a convergence result based on the Barles--Souganidis theorem \cite{BS}. In spite of its robustness, however, value iteration is a relatively inefficient technique to compute the numerical solution, and an acceleration strategy would be highly desirable. \\
From the very start of Dynamic Programming techniques \cite{B, H}, policy iteration (PI) has been recognized as a viable, usually faster alternative to value iteration in computing the fixed point of the Bellman operator. Among the wide literature on policy iteration, we quote here the pioneering theoretical analysis of Puterman and Brumelle \cite{PutBru}, which have shown that the linearization procedure underlying policy iteration is equivalent to a Newton-type iterative solver. More recently, the abstract setting of \cite{PutBru} has been adapted to computationally relevant cases \cite{SanRus}, proving superlinear (and, in some cases, quadratic) convergence of policy iteration. Moreover, we mention that an adaptation of policy iteration to large sparse problems has been proposed as ``modified policy iteration'' (MPI) in \cite{PutShi}, and has also become a classical tool.

In the present paper, we intend to study the construction and numerical validation of a SL scheme with PI/MPI sover for hybrid optimal control. To this end, we will recall the general algorithm, sketch some implementation details for the simple case of one-dimensional dynamics, and test the scheme on some numerical examples in dimension $d=1,2$.

The outline of the paper is the following. In Section \ref{costruz} we will review the main results about the Bellman equation characterizing the value function, and construct a semi-Lagrangian (SL) approximation for $V$ in the form of value iteration. In Section \ref{policy} we will improve the algorithm by a policy iteration technique. Finally, section \ref{tests} will present some numerical examples of approximation of the value function and construction of the optimal control.

\section{A Semi-Lagrangian scheme for hybrid control problems}\label{costruz}

First, we recall some basic analytical results about the value function \eqref{f_valore}. To this end, we start by making a precise set of assumptions on the problem.

\subsection{Basic assumptions and analytical framework}

In the product space $\R^d\times\I$, we consider sets (and in particular the sets $A, C$ and $D$) of the form
\begin{equation}\label{sets}
S = \{(x,q)\in\R^d\times\I: x\in S_i, q=i \},
\end{equation}
in which $S_i$ represents the subset of $S$ in which $q=i$. We assume that:

\begin{itemize}

\item[\bf (A1)] For each $q\in\I$, $A_q$, $C_q$, and $D_q$ are closed subsets of $\R^d$, and $D_q$ is bounded. $\partial A_q$ and $\partial C_q$ are $C^2$.

\item[\bf (A2)]  The function $f$ is bounded. Moreover, it is Lipschitz continuous in the state variable $x$ and uniformly continuous in the control variable $\alpha$.

\item[\bf (A3)] The map $g: A \times \mathcal{V}\to D$ is bounded and uniformly Lipschitz continuous with respect to $x$.

\item[\bf (A4)] $\partial A$ is a compact set, and for some $\gamma>0$, the following transversality condition:
\[
f(x,q,\alpha)\cdot \eta_{x,q} \leq -2\gamma
\]
holds for all $x\in \partial A_q$, and all $\alpha\in U$, where $\eta_{x,q}$ denotes the unit outward normal to $\partial A_q$ at $x$. We also assume similar transversality conditions on $\partial C$.

\item[\bf (A5)] We assume that, for all $i\in\I$, $d(A_i,C_i) \geq \beta>0$ and $d(A_i,D_i)\geq \beta >0$, where $d$ is the Euclidean distance.

\item[\bf (A6)] The control set $U$ is a compact metric space, and $\mathcal{V}$ is a finite discrete set.

\item[\bf (A7)] $\ell:\R^d\times \I\times U$ is a bounded and non-negative function, Lipschitz continuous w.r.t. the $x$ variable, and uniformly continuous w.r.t. the $\alpha$ variable.

\item[\bf (A8)] $c_A(x,q,w)$ and $c_C(x,q,x',q')$ are uniformly Lipschitz continuous in the variables $x$ and $x'$, and bounded with a strictly positive infimum. Moreover, for any $x$ and $q$, the function $c_C$ satisfies (for some $\Delta\ge 0$) the inequality
\[
c_C(x,q,x',q') < c_C(x,q,\bar x,\bar q) + c_C(\bar x,\bar q,x',q') - \Delta
\]

\end{itemize}

Via a suitable generalization of the Dynamic Programming Principle, it can be proved that the Bellman equation of the problem is in the form of a Quasi-Variational Inequality, and more precisely, once defined the Hamiltonian by
$$
 H(x,q,p) := \sup_{\alpha\in U}\big\{ -\ell(x,q,\alpha) - f(x,q,\alpha)\cdot p\big\}
$$
and the transition operators $\mathcal{M}$ and $\mathcal{N}$ by:
\begin{align*}
\mathcal{M} \phi(x,q) &:= \inf_{w\in \mathcal{V}} \Big\{\phi\big(g(x,q,w)\big)+c_A(x,q,w)\Big\} \qquad (x,q)\in A\\
\mathcal{N} \phi(x,q) &:= \inf_{(x',q')\in D} \big\{\phi(x',q')+c_C(x,q,x',q')\big\} \qquad (x,q)\in C
\end{align*}
we have the following

\begin{Theorem}[\cite{DhaRam}] Assume (A1)--(A8). Then, the function $V$ is the unique bounded and H\"older continuous viscosity solution of:
  \begin{equation}\label{hjb}
  \begin{cases}
  \lambda V(x,q) + H\big(x,q,D_x V(x,q)\big) = 0 & (x,q)\in(\R^d\times \I) \setminus (A\cup C) \\
  \max\Big\{V(x,q)-\mathcal{N} V(x,q),V(x,q) + H\big(x,q,D_x V(x,q)\big)\Big\} = 0 & (x,q)\in C \\
  V(x,q) - \mathcal{M} V(x,q) =0   & (x,q)\in A
  \end{cases}
  \end{equation}
\end{Theorem}

Note that uniqueness follows from a strong comparison principle, which also allows to use the Barles--Souganidis theorem \cite{BS} for proving convergence of stable and monotone schemes.

\subsection{Numerical approximation}

In order to set up a numerical approximation for \eqref{hjb}, we construct a discrete grid of nodes $(x_j,q)$ in the state space and fix the discretization parameters $\Delta x$ and $\Delta t$. In what follows, we will denote the discretization steps in compact form by $\delta:=(\Delta t, \Delta x)$ and the approximate value function by $V_\delta$.

Following \cite{FerZid}, we write the fixed point form of the scheme at $(x_i,q)$ as
\begin{equation}\label{scheme1}
v_i^{(q)} = V_\delta(x_i,q) =
\begin{cases}
\min \big\{ N V_\delta(x_i,q), \Sigma(x_i,q,V_\delta) \big\} & (x_i,q)\in C\\
M V_\delta(x_i,q) & (x_i,q)\in A\\
\Sigma(x_i,q,V_\delta) & \text{else}
\end{cases}
\end{equation}
in which $N$, $M$ and $\Sigma$ are consistent and monotone numerical approximations for respectively the operators $\mathcal{N}$, $\mathcal{M}$ and the Hamiltonian $H$. More compactly, \eqref{scheme1} could be written as
\[
V_\delta = T_\delta(V_\delta).
\]
We recall that, for $\lambda>0$, under the basic assumption which ensure continuity of the value function, the right-hand side of \eqref{scheme1} is a contraction \cite{FerZid} and can therefore be solved by fixed-point iteration, also known as {\it value iteration} (VI):
\begin{equation}\label{VI}
V_{\delta,j+1} = T_\delta(V_{\delta,j}).
\end{equation}
To define more explicitly the scheme, as well as to extend the approximate value function to all $x\in\R^d$ and $q\in\I$, we use an interpolation $\mathcal{I}$ constructed on the node values, and denote by $\mathcal{I}[V_\delta] (x,q)$ the interpolated value of $V_\delta$ computed at $(x,q)$. With this notation, a natural definition of the discrete jump operators $M$ and $N$ is given by
\begin{align}
 M V_\delta(x,q) &:= \min_{w\in\mathcal{V}}\Big\{\mathcal{I}[V_\delta]\big(g(x,q,w)\big) + c_A(x,q,w)\Big\}\label{M} \\
 N V_\delta(x,q) &:= \min_{(x',q')\in D}\Big\{\mathcal{I}[V_\delta](x',q') + c_C(x,q,x',q')\Big\} \label{N}
\end{align}
On the other hand, a standard semi-Lagrangian discretization of the Hamiltonian related to continuous control is provided (see \cite{FF}) by
\begin{equation}\label{scheme2}
\Sigma(x_i,q,V_\delta):= \min_{\alpha\in U}\Big\{\Delta t \> \ell(x_i,q,\alpha) + e^{-\lambda \Delta t}\>\mathcal{I}[V_\delta] \big(x_i+\Delta t \> f(x_i,q,\alpha), q\big)\Big\}.
\end{equation}
In the SL form, the value iteration \eqref{VI} might then be recast at a node $(x_i,q)$ as
\begin{equation}\label{scheme3}
v_{i,j+1}^{(q)} =
\begin{cases}
\displaystyle\min_{w\in\mathcal{V}}\Big\{\mathcal{I}[V_{\delta,j}]\big(g(x_i,q,w)\big) + c_A(x_i,q,w)\Big\} &(x_i,q)\in A \\
\displaystyle\min \Big\{ \min_{(x',q')\in D}\big\{\mathcal{I}[V_{\delta,j}](x',q') + c_C(x_i,q,x',q')\big\}, \Sigma(x_i,q,V_{\delta,j}) \Big\} & (x_i,q)\in C \\
\displaystyle\Sigma(x_i,q,V_{\delta,j}) & \text{else}
\end{cases}
\end{equation}
with $\Sigma$ given by \eqref{scheme2}, and $j$ denoting the iteration number.

Convergence of the scheme can be proved by using the arguments in \cite{FF, FerZid}) if the interpolation $\mathcal{I}$ is monotone (e.g., a $\P_1$ or $\Q_1$ finite element interpolation):

\begin{Theorem}[\cite{FerZid}] Assume (A1)--(A8). Assume in addition that $\lambda>0$, and that the interpolation $\mathcal{I}$ is monotone and invariant for the sum of constants. Then, $V_{\delta,j}\to V_\delta$ for $j\to\infty$. Moreover, the approximate solution $V_\delta$ converges to $V$ locally uniformly in $\R^d\times\I$ for $\Delta x,\Delta t\to 0$.
\end{Theorem}

\section{Policy iteration algorithm}\label{policy}

Following \cite{S}, we give now an even more explicit form of the scheme, which is the one applied to the one-dimensional examples of Sec.~\ref{tests}. Once we set up a 1-D space grid of evenly spaced nodes $x_1,\ldots,x_n$ with space step $\Delta x$, the discrete solution may be given the vector structure
$$
\bm v :=(\bm v ^{(1)},\bm v ^{(2)},\dots,\bm v ^{(m)}) \in \mathbb{R}^{n m}
$$
in which $\bm v ^{(q)}:=(v^{(q)}_1,\ldots,v^{(q)}_n)$ denotes the discretized value function associated to the $q$-th component of the state space. Within the vector $\bm v$, the element $v^{(k)}_i$ appears with the index $(k-1)n + i$.

Keeping the same notation for all vectors, $\bm\alpha \in U^{nm}$ will denote the vector of controls of the system, $\alpha_i^{(k)}$ being the value of the control at the space node $x_i$ while the $k$-th dynamics is active. We also define the vector $\bm s \in \I^{nm}$ representing the switching strategy, so that $s_i^{(k)}=l$ means that if the trajectory is in $x_i$ and the active dynamics is $k$, the system commutes from $k$ to $l$. Note that, in the numerical examples of Sec.~\ref{tests}, discontinuous jumps will always appear only on the discrete component of the state space, so that, for example, we have $x'=x$ and this data need not be kept in memory (we will use the term {\em switch} to denote a state transition of this kind).

In the general case, we would also need to keep memory of the arrival point of the jump and/or of the discrete control $w$ in the case of an autonomous jump. In general, the arrival point is not a grid point, so that we also need to perform an interpolation in \eqref{M}--\eqref{N}. Therefore, the details for the general case can be recovered by mixing the basic arguments used in what follows.

The endpoint of this construction is to put the problem in the standard form used in policy iteration,
\begin{equation}\label{eq:PI}
\min_{(\bm\alpha,\bm s) \in U^{nm}\times \I^{nm}}\big(B(\bm\alpha,\bm s)\bm v - \bm c(\bm \alpha,\bm s)\big)=0,
\end{equation}
with explicitly defined matrix $B$ and vector $\bm c$. Note that, in \eqref{eq:PI}, we have made clear the fact that a policy is composed of both a feedback control $\bm\alpha$ and a switching strategy $\bm s$.


Define now the matrices $D_A,D_C\in M_{nm}\big(\{0,1\}\big)$ as permutations of the array $\bm v$. These matrices represent changes in the state due to the switching strategy: $D_A$ corresponds to autonomous jumps and $D_C$ to controlled jumps. Note that, in our case, the elements of $D_A$ and $D_C$ will be in $\{0,1\}$, that there exists at most one nonzero element on each row, and that the two matrices cannot have a nonzero element in the same position.\\
In order to determine the positions of the nonzero elements ${d_A}_{a,b}(s) = 1$ in the matrix $D_A$, we apply the following rule. For all $(i,k) \in \{1,\ldots,n\}\times \I$, if the following conditions hold:
$$
\begin{cases}
(x_i,k) \in A \\
s_i^{(k)} \neq k & \text{(a switch occurs)}\\
(x_i,s_i^{(k)}) \in g(x_i,k,\mathcal{V}) & \text{(the switch is in the image of $g$),}
\end{cases}
$$
then,
$$
\begin{cases}
a=(k-1)n+i \\
b=(s_i^{(k)}-1)n+i. \\
\end{cases}
$$
Similarly, the nonzero elements of the matrix $D_C$, ${d_C}_{a,b}(s) = 1$, follow a slightly less strict rule. For all $\{1,\ldots,n\}\times \I$, if the following conditions hold:
$$
\begin{cases}
(x_i,k) \in C \\
s_i^{(k)} \neq k & \text{(a switch occurs),}
\end{cases}
$$
then,
$$
\begin{cases}
a=(k-1)n+i \\
b=(s_i^{(k)}-1)n+i. \\
\end{cases}
$$
Last, we define the matrix
$$
D(\bm s) := D_A(\bm s) + D_C(\bm s),
$$
which accounts for changes in the state related to the switching strategy, both autonomous and controlled.

We turn now to the continuous control part. First, we write $\Sigma$ in vector form as
$$
\bm\Sigma(\bm x,k,\bm v)= \min_{\bm\alpha^{(k)}\in U^n}\Big\{\Delta t \> \bm\ell(\bm x,k,\bm \alpha^{(k)}) + e^{-\lambda \Delta t}\mathcal{E}(\bm x,k,\bm\alpha^{(k)})\bm v^{(k)}\Big\},
$$
where the matrix $\mathcal{E}(\bm x,k,\bm\alpha^{(k)})  \in M_{n}(\mathbb{R})$ is defined so as to have
$$
\mathcal{E}(\bm x,k,\bm\alpha^{(k)})\bm v^{(k)} = \bm{\mathcal{I}}[V_\delta] \big(\bm x+\Delta t\> f(\bm x,k,\bm\alpha^{(k)}), k \big)
$$
and $\bm{\mathcal{I}}[V_\delta] \big(\bm x+\Delta t \> f(\bm x,k,\bm\alpha^{(k)}), k\big)$ and $\bm\ell(\bm x,k,\bm\alpha^{(k)})$ denote vectors which collect respectively all the values $\mathcal{I}[V_\delta] \big(x_i+\Delta t\> f(x_i,k,\alpha_i^{(k)}), k\big)$ and $\ell(x_i,k,\alpha_i^{(k)})$.

At internal points, using a monotone $\mathbb{P}_1$ interpolation for the values of $\bm v$ results in a convex combination of node values. On the boundary of the domain, the well-posedness of the problem requires either to have an invariance condition (which implies that $f(x_i,k,\alpha_i^{(k)})$ always points inwards) or to perform an autonomous jump or switch when the boundary is reached. Therefore, we should not care about defining a space reconstruction outside of the computational domain, although this could be accomplished by extrapolating the internal values.

The matrix $E(\alpha,s)\in M_{n m}(\mathbb{R})$ is then constructed in the block diagonal form:
$$
E(\bm\alpha,\bm s):=
\begin{pmatrix}
E^{(1)}(\bm\alpha^{(1)},\bm s^{(1)}) & 0 & \cdots & 0 \\
0 & E^{(2)}(\bm\alpha^{(2)},\bm s^{(2)}) & \ddots & \vdots\\
\vdots & \ddots & \ddots & 0\\
0 & \cdots & 0 & E^{(m)}(\bm\alpha^{(m)},\bm s^{(m)})
\end{pmatrix}
$$
Assuming for simplicity that we work at Courant numbers below the unity (although this is not necessary for the stability of SL schemes), each block $E^{(k)}(\bm\alpha^{(k)},\bm s^{(k)})\in M_n(\mathbb{R})$ is a sparse matrix with non-zero elements $e^{(k)}_{i,j}$ determined so as to implement a $\mathbb{P}_1$ space interpolation, in the following way: for every $(i,k) \in \{1,\ldots,n\}\times \I$, define
$$
h_{i,k}:=\frac{\Delta t}{\Delta x}  f(x_i,k,\alpha_i^{(k)})
$$
and
$$
\text{if}
\quad
\begin{cases}
s_i^{(k)}=k\\
h_{i,k}<0
\end{cases}
\quad
\text{then}
\quad
\begin{cases}
e^{(k)}_{i,i-1}(\alpha,s)=1+h_{i,k}\\
e^{(k)}_{i,i}(\alpha,s)=-h_{i,k}
\end{cases}
$$
else,
$$
\text{if}
\quad
\begin{cases}
s_i^{(k)}=k\\
h_{i,k}>0
\end{cases}
\quad
\text{then}
\quad
\begin{cases}
e^{(k)}_{i,i}(\alpha,s)=1-h_{i,k}\\
e^{(k)}_{i,i+1}(\alpha,s)=h_{i,k}
\end{cases}
$$
Note that, if a switching strategy $\bm z \in \I^{nm}$ doesn't perform any switch (i.e. $z_i^{(k)} = k$ for all $(i,k) \in N\times\I $), by definition of the matrix $E(\alpha,s)$ we would obtain, for all $k \in \I$,
$$
E^{(k)}(\bm\alpha^{(k)},\bm z^{(k)})\bm v^{(k)}=\mathcal{E}(\bm x, k, \bm\alpha^{(k)})\bm v^{(k)},
$$
whereas, in the general case, if a switch occurs at $x_i$, then the corresponding element of the matrix $E^{(k)}$ is zero.
Finally, we define the vector $\bm c(\bm\alpha,\bm s) \in \mathbb{R}^{n m}$ with a block structure of the form :
 $$
\bm c(\bm\alpha,\bm s)=\big(\bm c^{(1)}(\bm\alpha^{(1)},\bm s^{(1)}),\bm c^{(2)}(\bm\alpha^{(2)},\bm s^{(2)}),\dots,\bm c^{(m)}(\bm\alpha^{(m)},\bm s^{(m)})\big)
$$
with $\bm c^{(k)}(\bm\alpha^{(k)},\bm s^{(k)}) \in \mathbb{R}^n$ such that, for every $(i,k)$ in $\{1,\ldots,n\}\times \I$,
$$
c^{(k)}_i(\alpha_i^{(k)},s_i^{(k)})=
\begin{cases}
-\Delta t\>\ell(x_i,k,\alpha_i^{(k)})& s_i^{(k)}=k\\
-\xi(k,s_i^{(k)})& s_i^{(k)}\neq k\\
\end{cases}
$$
where $\xi(k,l)$ denotes the switching cost ($c_A$ or $c_C$) from dynamics $k$ to $l$.
 
With these notations, we can write the SL scheme \eqref{scheme1} in vector form as
\begin{equation}\label{scheme4}
\bm v = \min_{(\bm\alpha,\bm s) \in U^{nm}\times \I^{nm}}\Big\{\big[D(\bm s)+e^{-\lambda \Delta t}E(\bm\alpha,\bm s)\big]\bm v-c(\bm\alpha,\bm s)\Big\},
\end{equation}
or, defined the matrix
$$
B(\bm\alpha,\bm s):=-I_{nm}+D(\bm s)+e^{-\lambda \Delta t}E(\bm\alpha,\bm s),
$$
as
 $$
\min_{(\bm\alpha,\bm s) \in U^{nm}\times \I^{nm}}\big(B(\bm\alpha,\bm s)\bm v-\bm c(\bm\alpha,\bm s)\big)=0.
 $$
Once we have reformulated the Semi-Lagrangian scheme for the hybrid control problem in the standard form, we can solve it using Algorithm \ref{how}. The only difference with a standard PI algorithm is to include the switching strategy in the control policy.

\begin{algorithm}
\begin{algorithmic}
\State $j\gets 0$
\State STOP $\gets$ FALSE
\State $\bm\alpha_0 \in U^{nm}$
\State $\bm s_0 \in \I^{nm}$
\While {STOP = FALSE}
\If{[stopping criterion satisfied]}
\State STOP $\gets$ TRUE
\Else
\State $\bm v_{j}\gets \bm w$ solution of $B(\bm\alpha_j, \bm s_j)\bm w = \bm c(\bm\alpha_j, \bm s_j)$ \hspace\fill (\emph{Policy Evaluation})
\State $(\bm\alpha_{j+1},\bm s_{j+1})\gets \displaystyle\argmin_{(\bm a,\bm\sigma) \in U^{nm}\times \I^{nm}}\big(B(\bm a,\bm\sigma)\bm v_j - \bm c(\bm a,\bm\sigma)\big)$ \hspace\fill (\emph{Policy Improvement})
\State $j\gets j+1$
\EndIf
\EndWhile
\end{algorithmic}
\caption{Policy Iteration, 1D matrix form} \label{how}
\end{algorithm}

We remark that some theoretical result obtained in the ``classical'' setting is also true in the hybrid setting. In particular, convergence might still be obtained by monotonicity (see, e.g., \cite{FF}) with minor changes in the proof, since the right-hand side of \eqref{scheme3} is still in the form of a minimum:

\begin{Theorem}
Let $\bm v$ be the solution of \eqref{scheme4}, and $\bm v_j$ be defined by Algorithm \ref{how}. If
\[
\min_{\bm\alpha,\bm s}\big(B(\bm\alpha,\bm s)\bm v_0 - \bm c(\bm\alpha,\bm s)\big)\le 0,
\]
then the sequence $\bm v_j$ is monotone decreasing, and $\bm v_j\to \bm v$.
\end{Theorem}

\subsection{Modified policy iteration}

A different iterative solver for the numerical scheme has been first proposed and analysed in \cite{PutShi}, and is known as {\it modified policy iteration}. It consists in performing the minimization in \eqref{scheme3} only once every $N_{it}$ iterations. In other terms, the policy evaluation step is replaced by $N_{it}$ iterations of linear advection (in which, however, the transport may occur among different components of the state space). For $N_{it}=1$ we obtain the value iteration, whereas for $N_{it}\to\infty$ the transport steps converge to an exact policy evaluation, and the algorithm coincides with the previous ``exact'' PI algorithm.

The pseudo-code in Algorithm \ref{inexact} shows the MPI algorithm in one-dimensional matrix form, for a comparison with the exact algorithm (Algorithm \ref{how}).

\begin{algorithm}
\begin{algorithmic}
\State $j\gets 0$
\State STOP $\gets$ FALSE
\State $\bm\alpha_0 \in U^{nm}$
\State $\bm s_0 \in \I^{nm}$
\While {STOP = FALSE}
\If{[stopping criterion satisfied]}
\State STOP $\gets$ TRUE
\Else
\If{$j=0 \>\> (\mathrm{mod}\> N_{it})$}
\State $(\bm\alpha_{j+1},\bm s_{j+1})\gets \displaystyle\argmin_{(\bm a,\bm\sigma) \in U^{nm}\times \I^{nm}}\big(B(\bm a,\bm\sigma)\bm v_j - \bm c(\bm a,\bm\sigma)\big)$ \hspace\fill (\emph{Policy Improvement})
\Else
\State $(\bm\alpha_{j+1},\bm s_{j+1})\gets (\bm\alpha_{j},\bm s_{j})$
\EndIf
\State $\bm v_{j+1}\gets \left[D(\bm s_{j+1})+e^{-\lambda \Delta t}E(\bm\alpha_{j+1},\bm s_{j+1})\right]\bm v_{j} - \bm c(\bm\alpha_{j+1},\bm s_{j+1})$ \hspace\fill (\emph{Inexact Policy\\ \hspace{12cm} Evaluation})
\State $j\gets j+1$
\EndIf
\EndWhile
\end{algorithmic}
\caption{Modified Policy Iteration, 1D matrix form} \label{inexact}
\end{algorithm}

Note that, in the numerical test section, the MPI algorithm has been applied to the two-dimensional examples. Although the formulation in dimension $d=2$ could be accomplished by a suitable redefinition of the vectors and matrices, in practice the MPI algorithm does not need such a formalism.

Concerning convergence, the hybrid case can again be treated with the same theoretical tools of the original proof in \cite{PutShi}, which relies on the monotonicity of the (discretized) Bellman operator, as well as on giving an upper and a lower bound on the sequence $\bm v_j$ by means of two converging sequences (one of which generated by value iteration). More precisely, the sequence considered in the convergence proof for the MPI is the sequence of approximations obtained after each policy improvement. In our notation, this is the subsequence $\bm v_l$ corresponding to $j=lN_{it}+1$. We have therefore the following

\begin{Theorem}
Let $\bm v$ be the solution of \eqref{scheme4}, and $\bm v_j$ be defined by Algorithm \ref{inexact}. If
\[
\min_{\bm\alpha,\bm s}\big(B(\bm\alpha,\bm s)\bm v_0 - \bm c(\bm\alpha,\bm s)\big)\le 0,
\]
then, for any $N_{it}\ge 1$, the subsequence $\bm v_l$ obtained for $j=lN_{it}+1$ is monotone decreasing, and $\bm v_l\to \bm v$ for $l\to\infty$.
\end{Theorem}

\section{Numerical tests}\label{tests}

We give in this section some numerical examples in one and two space dimensions, comparing the performances of Value and Policy Iteration -- exact PI algorithm in one dimension, and MPI in two dimensions. The comparison shows a substantial improvement in the convergence of the solver for the exact PI algorithm, whereas the MPI performs roughly the same number of iterations as the VI. Here, the bottleneck is apparently the contraction coefficient of the Bellman operator. Nevertheless, the MPI allows to avoid the minimization step in a large majority of the iterates, thus reducing the CPU time. Note that in both two-dimensional examples the control appears only as a switching strategy, and the complexity of policy evaluation steps is reduced by a factor $1/m$. For more complex control actions the improvement in computing time would be even greater.

\subsection{Stabilization of an unstable system}
We now apply this technique to a stabilization problem: we consider a system with two dynamics: one ``strong and expensive'' and the other ``weak and cheap''. Only the former is able to keep the state of the system within the given set over time.

The state equation $\dot X(t) = f\big(X(t), Q(t), \alpha(t)\big)$ is defined by
$$
f(x,q,\alpha)=
\begin{cases}
x+ d_1 \alpha & q=1 \\
x+ d_2 \alpha & q=2
\end{cases}
$$
where $d_1 <d_2$ and $-1 \leq \alpha(t) \leq 1$ for every $t$ in $[0,+\infty)$. Switching is mandatory only when the dynamics $q=1$ is active and $|X(t)|=1$, which implies that the state of the system belongs to the interval $[-1, 1]$ for all $t$ in $[0,+\infty)$.

Here and in what follows, $c_{i,j}$ denotes a constant switching cost from $q=i$ to $q=j$, and the cost functional is defined as 
$$
\ell(x,q,\alpha)=
\begin{cases}
x^2+c_1\alpha^2 & q=1 \\
x^2+c_2\alpha^2 & q=2
\end{cases}
$$
The values assigned to all the parameters are summed up in Table \ref{tab:ws_par}, whereas Table \ref{tab:ws} reports the number of iterations required for given stopping tolerances. In the first three examples, the stopping criterion reads
$$
\|\bm v_j - \bm v_{j-1}\|_\infty < \epsilon.
$$
Note that, according to Table \ref{tab:ws}, squaring the tolerance makes the number of iteration $N_P$ of the PI algorithm increase linearly, which indicates roughly quadratic convergence, while the number $N_V$ for VI has a geometric behaviour as expected.

\begin{table}
\centering
\begin{tabular}{|c|c|c|c|c|c|c|c|}
\hline
$d_1$ & $d_2$ & $c_{1,2}$ & $c_{2,1}$ & $c_1$ & $c_2$ & $\lambda$  & $t_f$ \\
\hline \hline
0.5 & 2 & 0.2 & 0 & 0.25 & 4 & 1& 20\\
\hline
\end{tabular}
\caption{Choice of parameters, weak-strong test}\label{tab:ws_par}
\bigskip
\begin{tabular}{|c|c|c|}
\hline
$\epsilon$ & $N_V$ & $N_P$ \\
\hline \hline
$10^{-3}$ & 456 &  8 \\
\hline
$10^{-6}$ & 1147 &  10 \\
\hline
$10^{-12}$ & 2786 &  12 \\
\hline
\end{tabular}
\caption{Number of iterations (VI and PI) for a given tolerance $\epsilon$, weak-strong test}
\label{tab:ws}
\end{table}

\begin{figure}
    \centering
    \includegraphics[width=\textwidth]{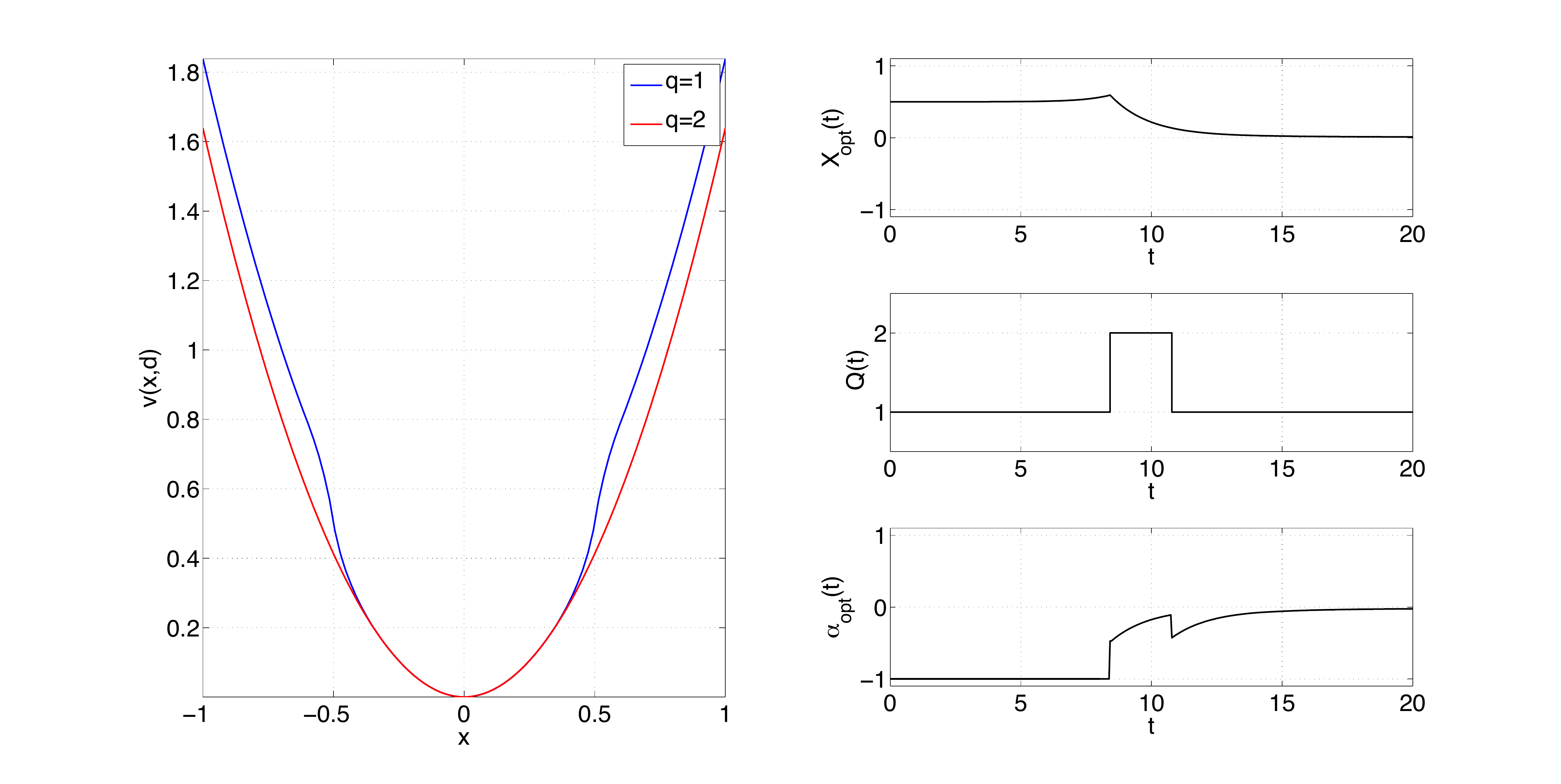}
    \caption{Value function, trajectory and optimal control, weak-strong test}
    \label{fig:ws}
\end{figure}

Figure \ref{fig:ws} shows the optimal strategy obtained for $\Delta t=0.0067$, $\Delta x=\Delta t||f||_{\infty}$, $\big(X(0),Q(0)\big)=(0.5,1)$, and $t\in[0, t_f]$. This strategy consists in using the unstable dynamics $q=1$ as long as the state belongs to the interval $[-\bar x, \bar x]$ (where the value of $\bar x \in [-1,1]$ depends on the given data). On the other hand, as soon as $|X(t)|>\bar x$, the optimal choice is to switch from $q=1$ to $q=2$ in order to stabilize the system and force it back towards the origin, then switch again to the first dynamics which can be used at a lower cost.

\subsection{Three-gear vehicle}

In this test, we consider the optimal control of a vehicle equipped with a three-gear engine, focusing on the acceleration strategy and the commutation between gears. Physical parameters correspond to the italian scooter Piaggio Vespa 50 Special.

The state equation for the speed of the vehicle is defined, for each gear $q \in \{1,2,3\}$, by
$$
f(x,q,\alpha):=\frac{1}{m}\left(\frac{T(\beta_q x)}{r \rho_q} \alpha-c_\text{d} x^2\right)
$$
where $T(\omega):=\tau \left(\frac{\omega}{\nu}-\left(\frac{\omega}{\nu}\right)^3\right)$ is the power band of the engine, $\tau$ and $\nu$ are respectively its maximum torque and r.p.m., $\beta_q:=\frac{60}{r \pi \rho_q}$ is a conversion coefficient with
$$
\rho_q:=\frac{\text{transmission shaft r.p.m.}}{\text{crankshaft r.p.m.}},
$$
$r$ is the radius of the wheel and $c_d$ the drag coefficient.
The control $\alpha(t) \in [0,1]$ represents the fraction of maximum torque used and the running cost is a linear combination of $x$ and $\alpha$:
$$
\ell(x,\alpha) = -c_x x+c_\alpha \alpha
$$
$c_x$ and $c_\alpha$ are positive weights.
Last, we define $c_{i,j}$ as the switching cost from from $q=i$ to $q=j$.

The numerical results are obtained by assigning realistic values (Table \ref{tab:parameter}) to the parameters defining the dynamics. The number of iterations is shown in Table \ref{tab:vespa} for various tolerances. The constant number of iterations obtained by PI might be due to the fact that optimal solutions (seem to) work with increasing values of $q$, this possibly meaning some sort of ``causality'' in the propagation of the value function.

\begin{figure}
    \centering
    \includegraphics[width=.6\textwidth]{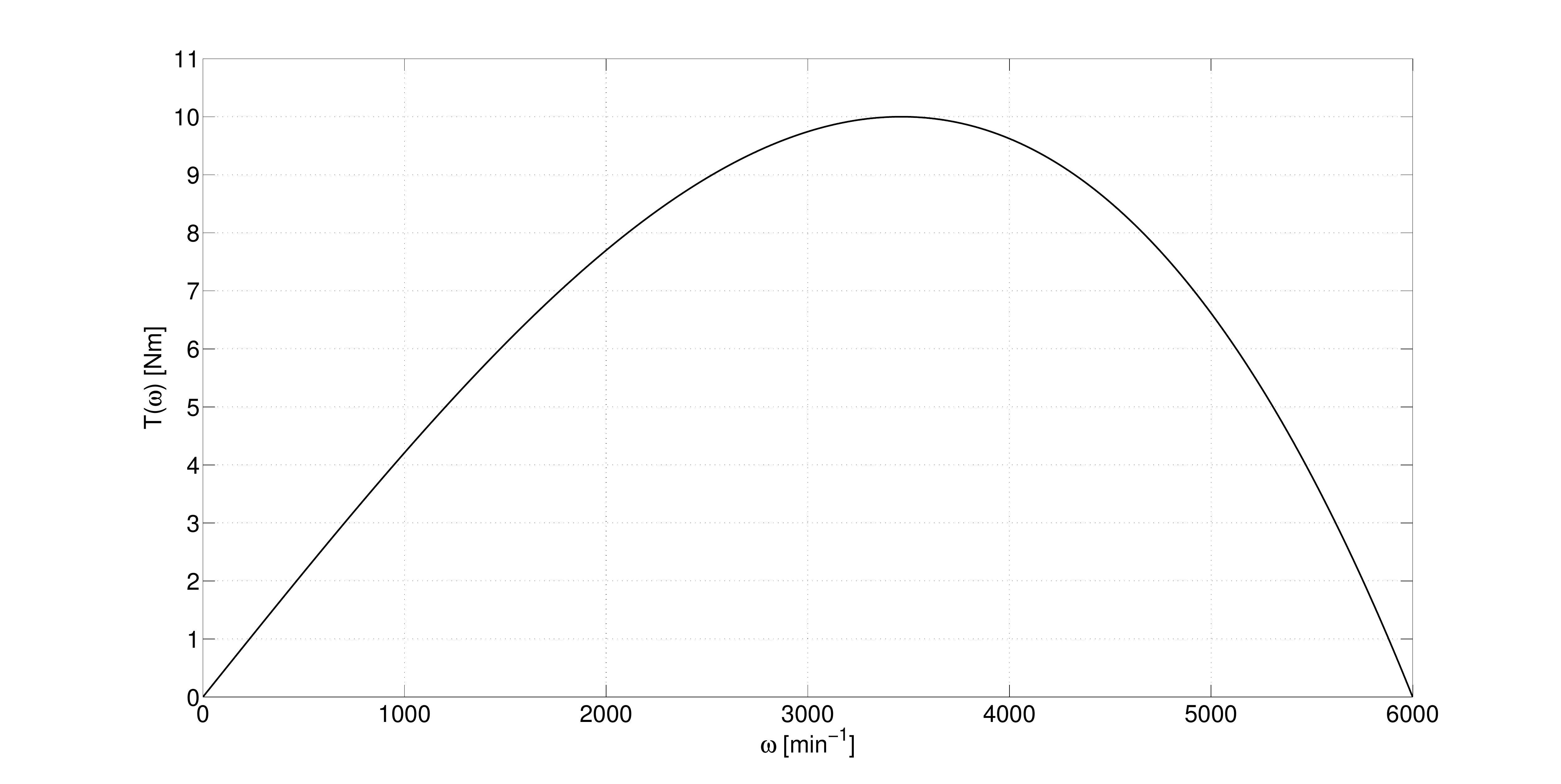}
    \caption{Power band of the engine.}
    \label{fig:tor}
    \centering
    \includegraphics[width=\textwidth]{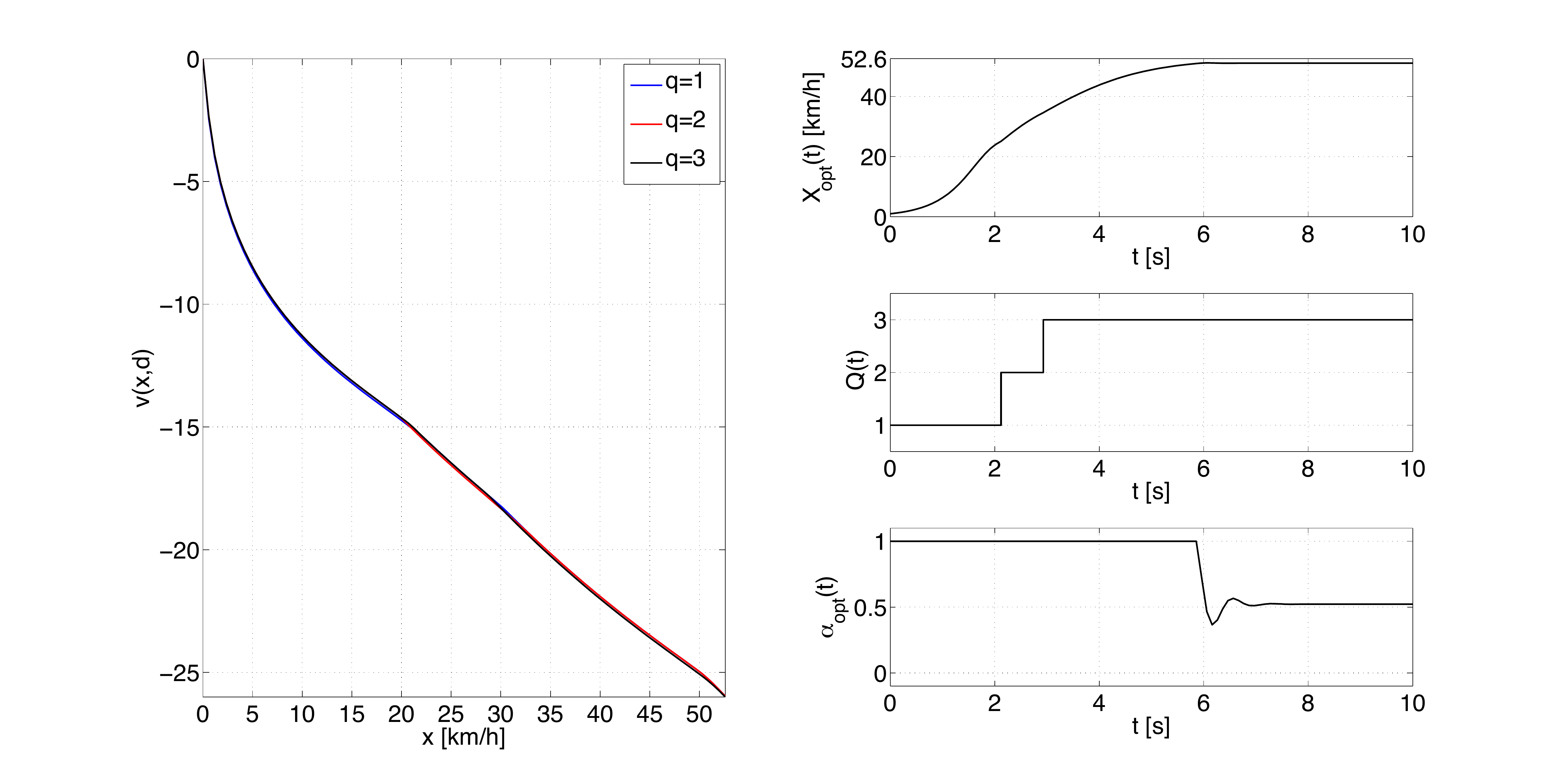}
    \caption{Value functions, trajectory and optimal control, three-gear vehicle, first case.}
    \label{fig:vespa}
    \centering
    \includegraphics[width=\textwidth]{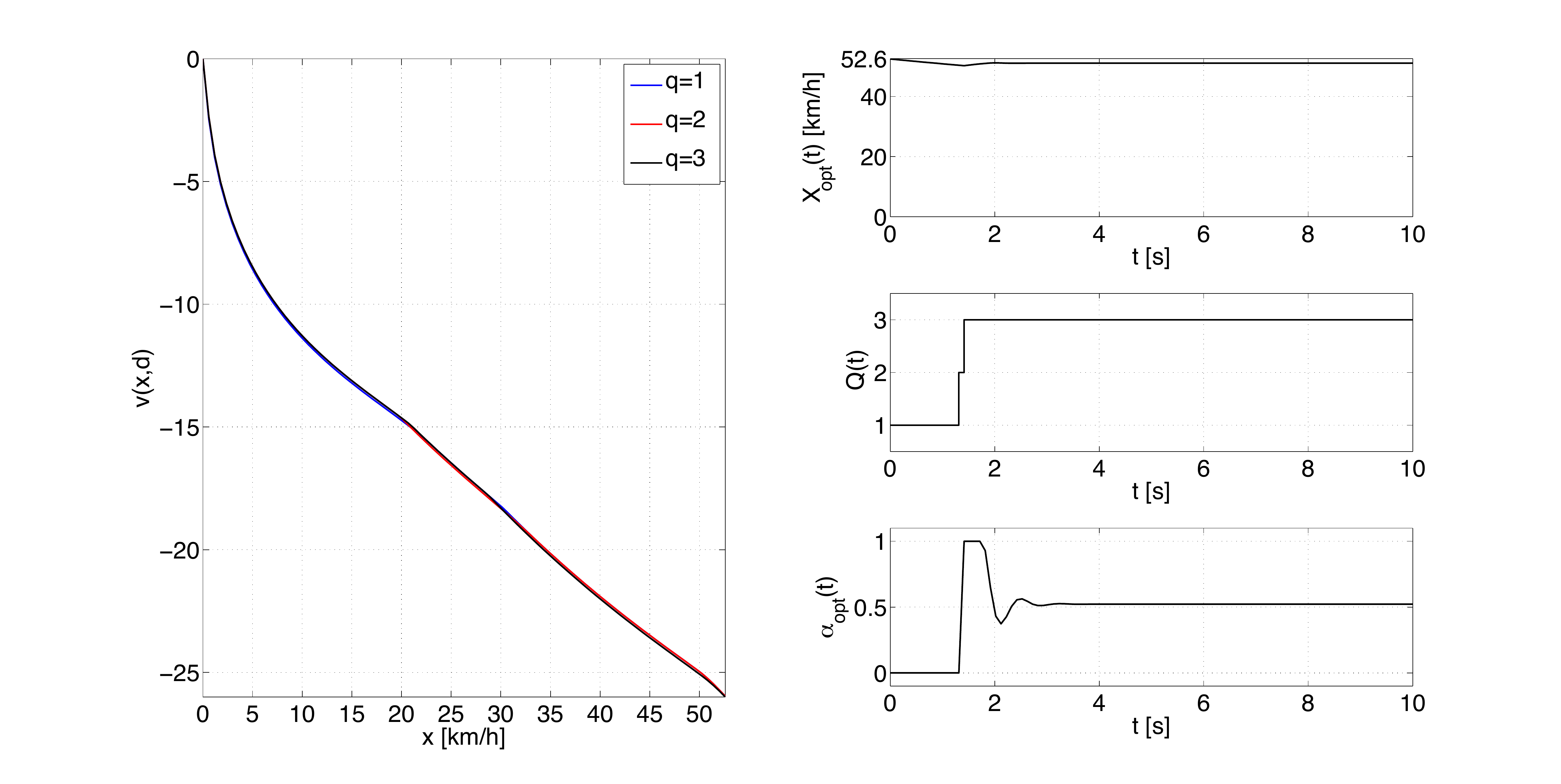}
    \caption{Value functions, trajectory and optimal control, three-gear vehicle, second case.}
    \label{fig:vespa2}
\end{figure}

\begin{table}
\centering
\begin{tabular}{|c|c|c|c|c|c|c|c|}
\hline
$m$ & $\rho_1$&$\rho_2$&$\rho_3$ & $r$ & $c_d$ & $\tau$ & $\nu$  \\
\hline \hline
$140$ [kg] & $0.06$&$0.09$&$0.12$ & $0.2$ [m] & $0.3$ & $ 10$ [Nm]& $6\cdot 10^3$ [min$^{-1}$] \\
\hline
\end{tabular}
\newline
\vspace*{0.5 cm}
\newline
\begin{tabular}{|c|c|c|c|c|}
\hline
$c_\alpha$ & $c_x$&$c_{i,j}$&$\lambda$& $t_f$ \\
\hline \hline
1 & 0.5 & $\begin{cases}0.1 & i \neq j\\ 0 & i = j\end{cases}$& 1 & $10$ [s] \\
\hline
\end{tabular}
\caption{Choice of parameters, three-gear vehicle test}\label{tab:parameter}
\bigskip
\begin{tabular}{|c|c|c|c|}
\hline
$\epsilon$ & $N_V$ & $N_P$ \\
\hline \hline
$10^{-3}$ & $337$ & $6$  \\
\hline
$10^{-6}$ & $586$ & $6$  \\
\hline
$10^{-12}$ & $1084$ & $6$   \\
\hline
\end{tabular}
\caption{Number of iterations (VI and PI) for a given tolerance $\epsilon$, three-gear vehicle test.}
\label{tab:vespa}
\end{table}
Figure \ref{fig:tor} shows the power band corresponding to our choice of $\tau$ and $\nu$.
Figure \ref{fig:vespa} shows the optimal solution obtained with $\Delta t=0.027$ [s], $\Delta x=\Delta t||f||_{\infty}$, $(x,q)=(0.28 \text{ [m/s]},1)$ and $t\in[0,t_f]$. The optimal strategy is to reach the highest gear as fast as possible and then stabilize at a value $\alpha\approx 0.5$. A different scenario is shown in Fig. \ref{fig:vespa2}, in which we set the initial state to $(x,q)=(14.58\text{ [m/s]},1)$. Here, the control lets the vehicle slow down, then switches to the third gear in order to replicate the previous strategy.

\subsection{Bang--Bang control of a chemotherapy model}

In this test, we consider the control of a two-compartment model of tumor growth. For this model, and cost functionals of the kind we will consider below, optimal controls are known to be of bang--bang type (see \cite{LedSch}). In this case, we can recast the problem in hybrid form, by considering an evolution in lack of chemotherapy ($Q=1$):
\begin{equation}\label{tum_model_1}
\begin{cases}
\dot X_1(t) = -a_1 X_1(t) + 2a_2X_2(t) \\
\dot X_2(t) = a_1 X_1(t) - a_2 X_2(t)
\end{cases}
\end{equation}
and a different evolution at full-dose chemotherapy ($Q=2$):
\begin{equation}\label{tum_model_2}
\begin{cases}
\dot X_1(t) = -a_1 X_1(t) \\
\dot X_2(t) = a_1 X_1(t) - a_2 X_2(t).
\end{cases}
\end{equation}
Here, the two compartments represent the number of cells at different stages of their lives, and the chemotherapy acts by preventing the generation of new tumor cells in the first compartment by inhibiting the mitosis of cells in the second compartment.

The cost functional is defined as
\begin{equation}\label{cost_chemo}
J(x,q;\theta) = \int_0^\infty \big(r_1 \dot X_1(t) + r_2 \dot X_2(t) + Q(t) - 1\big) e^{-\lambda t} dt,
\end{equation}
in which $\dot X_1(t)$ and $\dot X_2(t)$ are given by \eqref{tum_model_1}--\eqref{tum_model_2} for respectively $Q(t)=1$ and $Q(t)=2$, and we have to minimize a combination between the growth of the tumor mass and the toxic effect of the drug on healthy cells (note that this latter term appears only when $Q(t)=2$). Due to the geometric properties of the problem, Zeno executions cannot occur, and we can avoid to introduce a switching cost, which would have no practical meaning. Setting the switching cost to zero also causes the two value functions to coincide, i.e., $V(x,1)\equiv V(x,2)$, and in this case a switch can occur at $t=0^+$. While the general theory usually rules out this situation, no particular problems arise in this specific case.

\begin{figure}
    \centering
    \includegraphics[width=.49\textwidth]{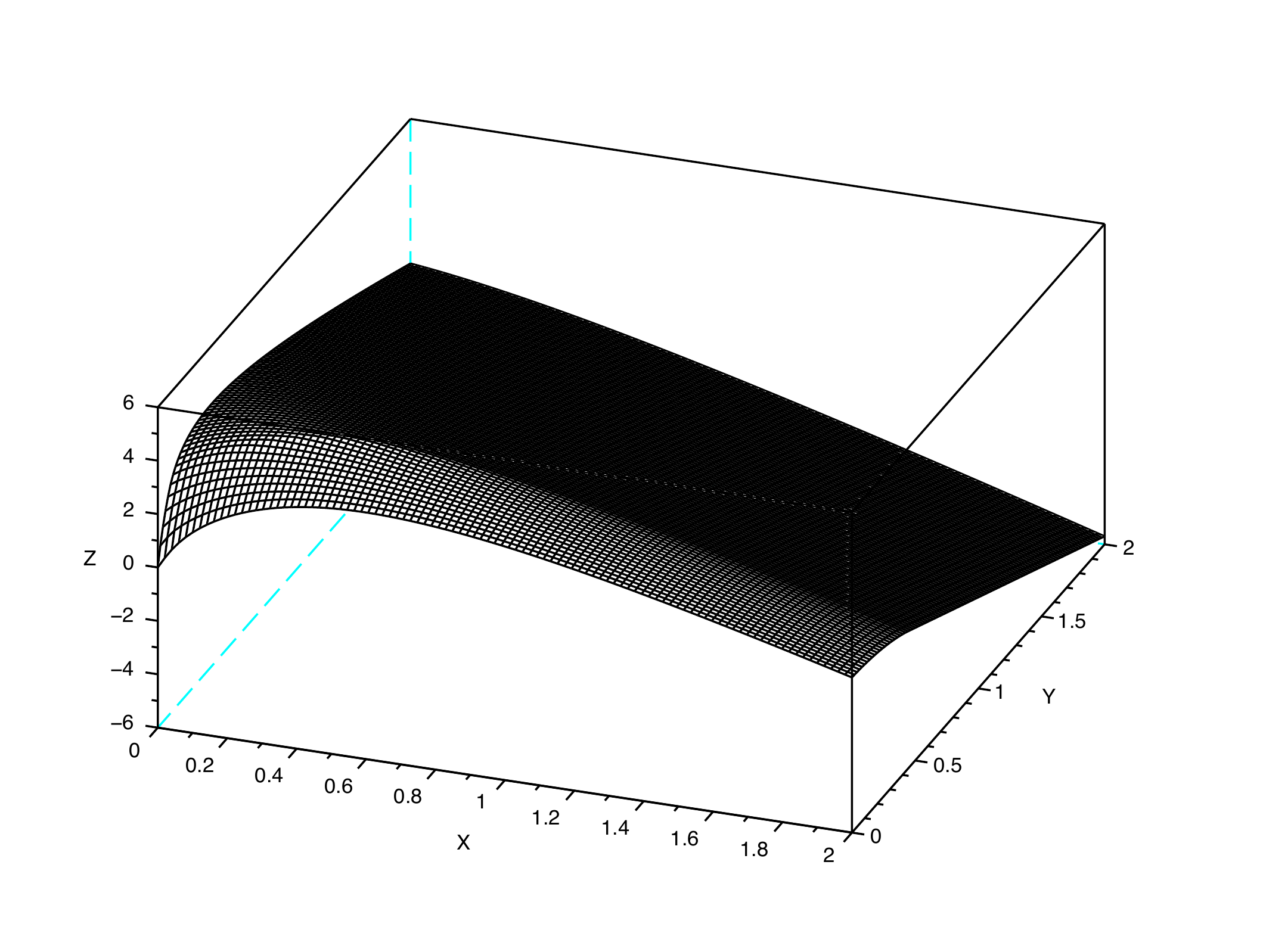}
    \includegraphics[width=.49\textwidth]{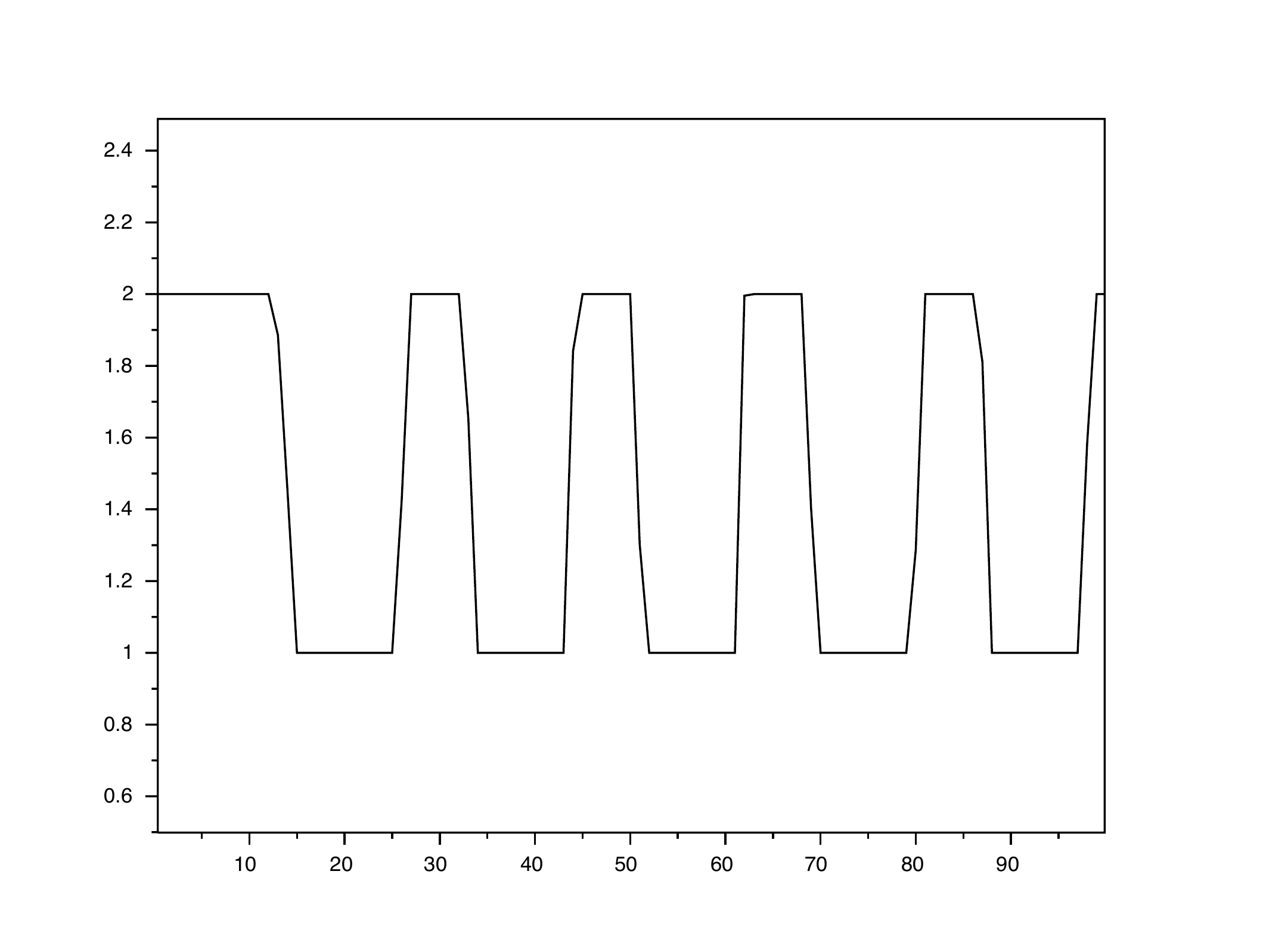}
    \caption{Value function and optimal switching for the chemotherapy test.}
    \label{fig:chemo1}
    \includegraphics[width=.49\textwidth]{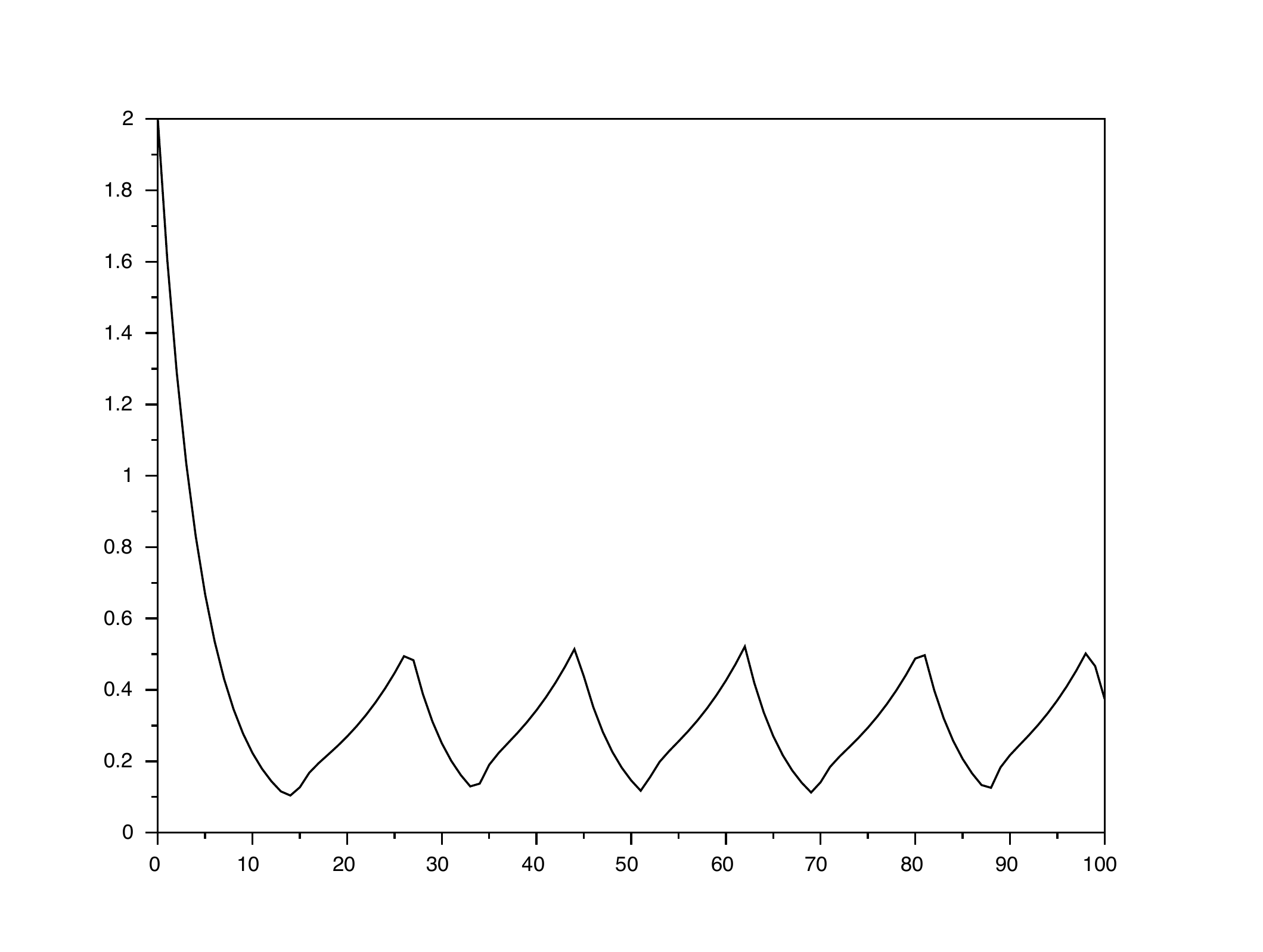}
    \includegraphics[width=.49\textwidth]{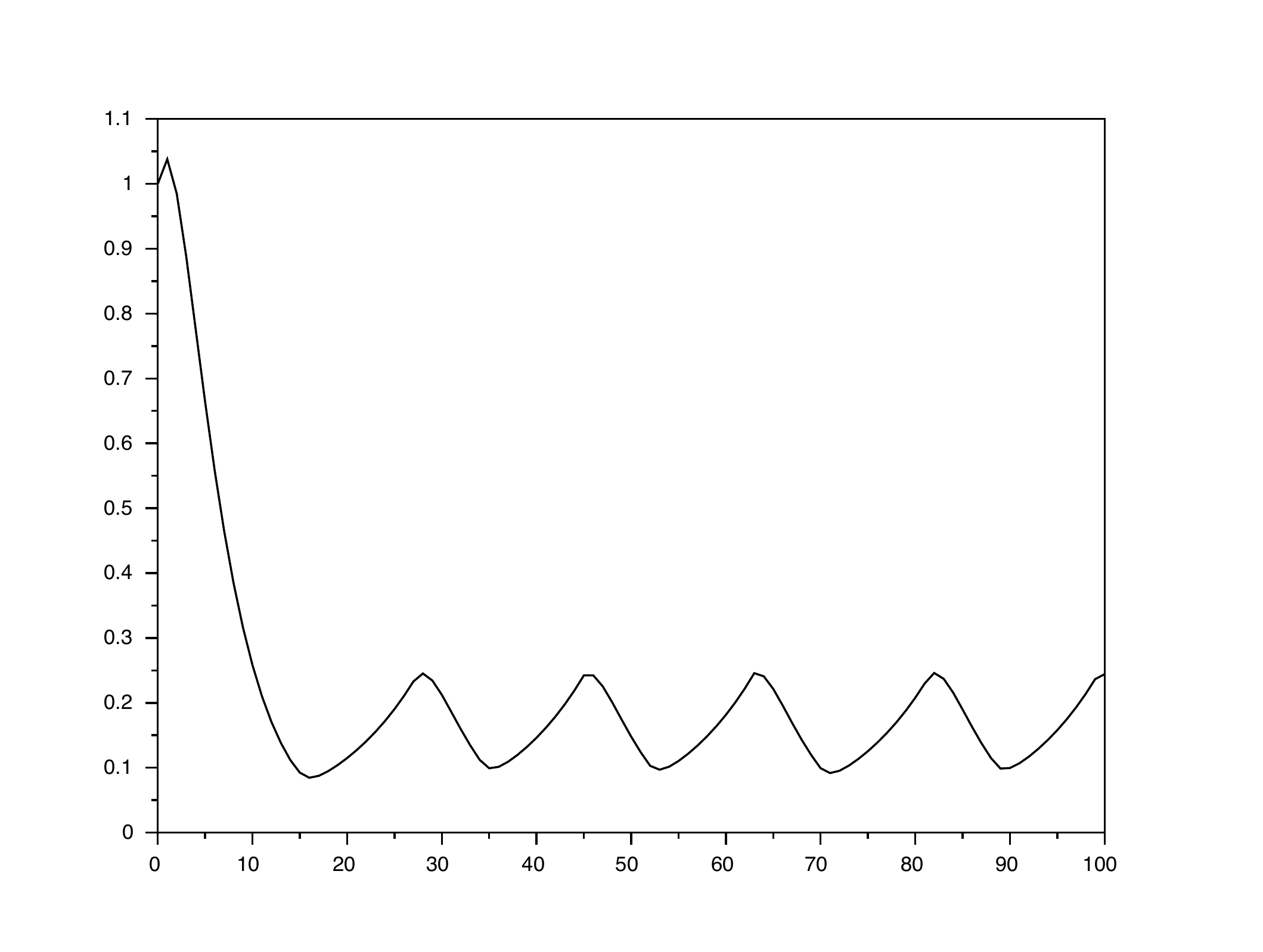}
    \caption{Trajectories $X_1(t)$ and $X_2(t)$ for the chemotherapy test.}
    \label{fig:chemo2}
    \includegraphics[width=.7\textwidth]{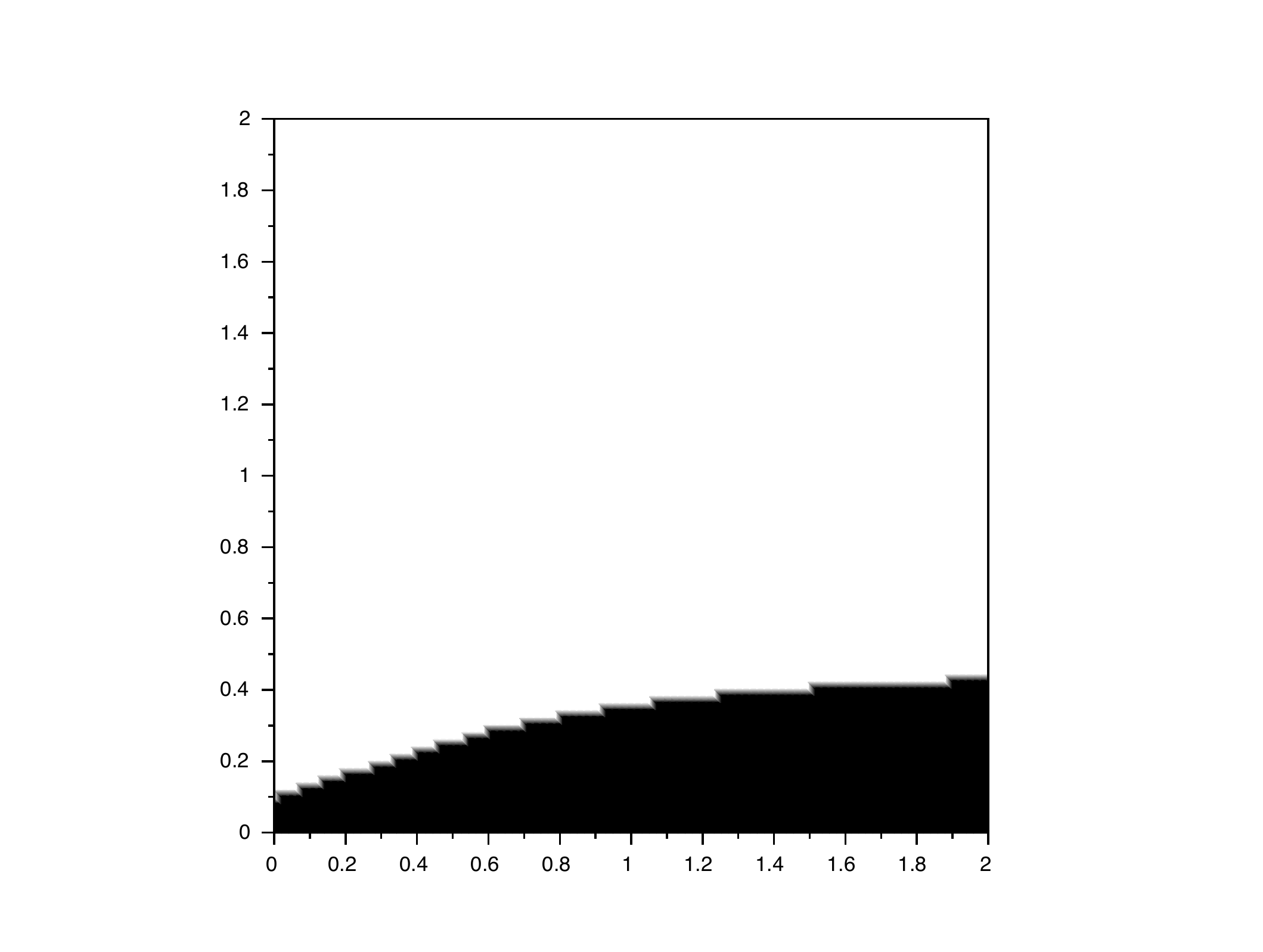}
    \caption{Optimal switching map for the chemotherapy test.}
    \label{fig:chemo3}
\end{figure}

The values of the parameters are assigned as in Table \ref{tab:parameter2}, according to the current literature (see \cite{LedSch}). Figg. \ref{fig:chemo1}--\ref{fig:chemo3} show the value function(s) of the problem, the optimal switching with respect to time and space and a sample trajectory starting from the initial state $(x_1,x_2) = (2,1)$. The value function has been computed with a $100\times 100$ grid on the domain $[0,2]^2$, and $\Delta t=0.1$. Note that there exists a clear discontinuity for the gradient of the value function, which corresponds to the switching curve in Fig. \ref{fig:chemo3}, which separates the black region, in which the optimal solution is $Q(t)=1$, from the white region, in which the optimal solution is $Q(t)=2$. The approximate optimal control shows a limit cycle in which a quasi-periodic switching between the two dynamics takes place.
\begin{table}
\centering
\begin{tabular}{|c|c|c|c|c|c|c|}
\hline
$a_1$ & $a_2$ & $r_1$ & $r_2$ & $\lambda$ & $x_1$ & $x_2$ \\
\hline \hline
$0.197$ & $0.356$ & $6.94$ & $3.94$ & $0.1$ & $2$ & $1$ \\
\hline
\end{tabular}
\caption{Choice of parameters, chemotherapy test}\label{tab:parameter2}
\bigskip
\begin{tabular}{|c|c|c|c|}
\hline
$\epsilon$ & $N_V$ & $N_P$ \\
\hline \hline
$10^{-3}$ & $192$ & $192$  \\
\hline
$10^{-6}$ & $528$ & $526$  \\
\hline
\end{tabular}
\caption{Number of iterations for a given tolerance $\epsilon$, chemotherapy test.}
\label{tab:chemo}
\end{table}

Table \ref{tab:chemo} compares the two (VI and MPI) numerical solvers. Here and in the following test, the MPI algorithm has been implemented with $N_{it}=10$, and an initial block of 10 value iterations has been performed at the very start in order to provide a better initial guess.
As remarked above, the Modified Policy Iteration algorithm performs essentially the same number of iterations than the value iteration algorithm, but at a lower cost.

\subsection{DC/AC inverter}

\begin{figure}
    \centering
    \includegraphics[width=.5\textwidth]{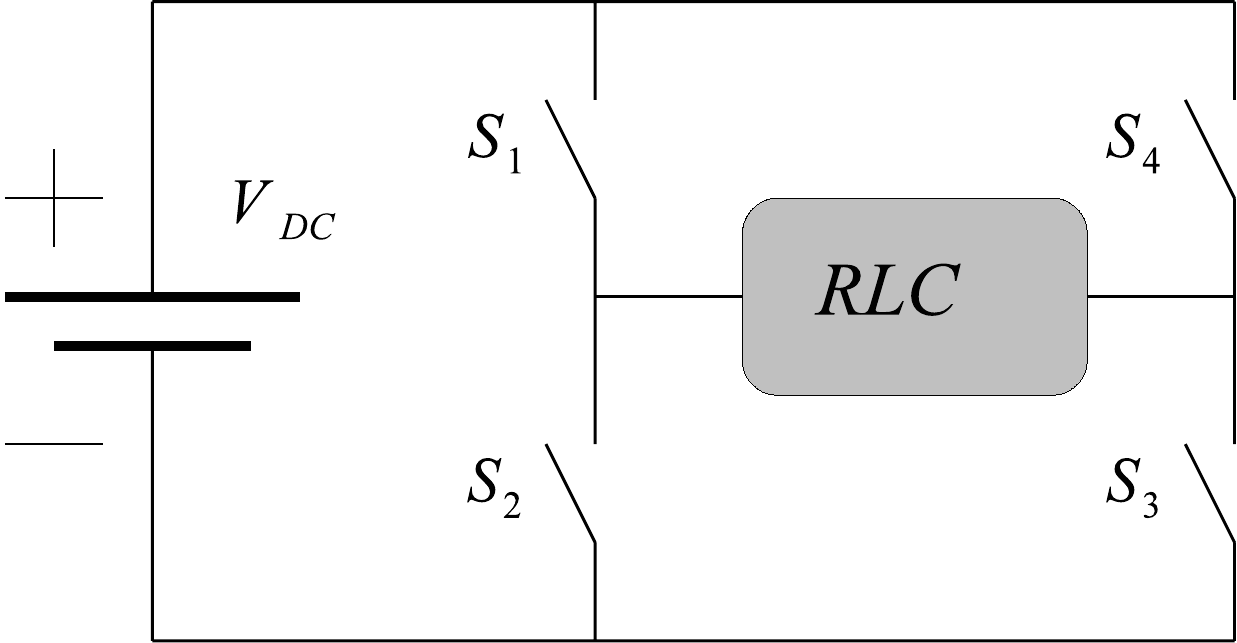}
    \caption{Abstract structure of the single-phase DC/AC inverter.}
    \label{fig:inverter_schema}
\end{figure}

The last test presents a single-phase DC/AC inverter, whose conceptual structure is sketched in Fig. \ref{fig:inverter_schema}.

\begin{figure}
    \centering
    \includegraphics[width=.49\textwidth]{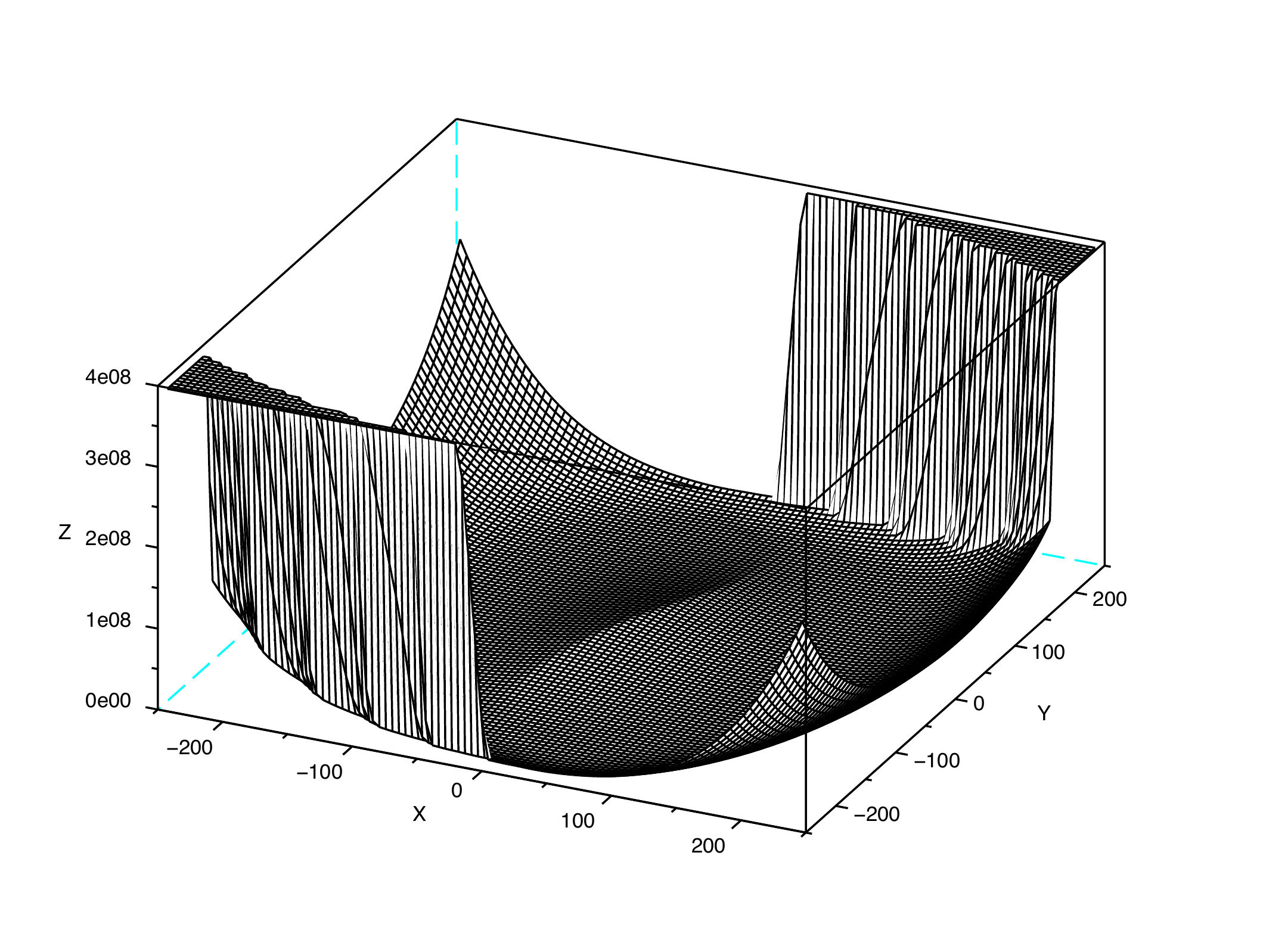}
    \includegraphics[width=.49\textwidth]{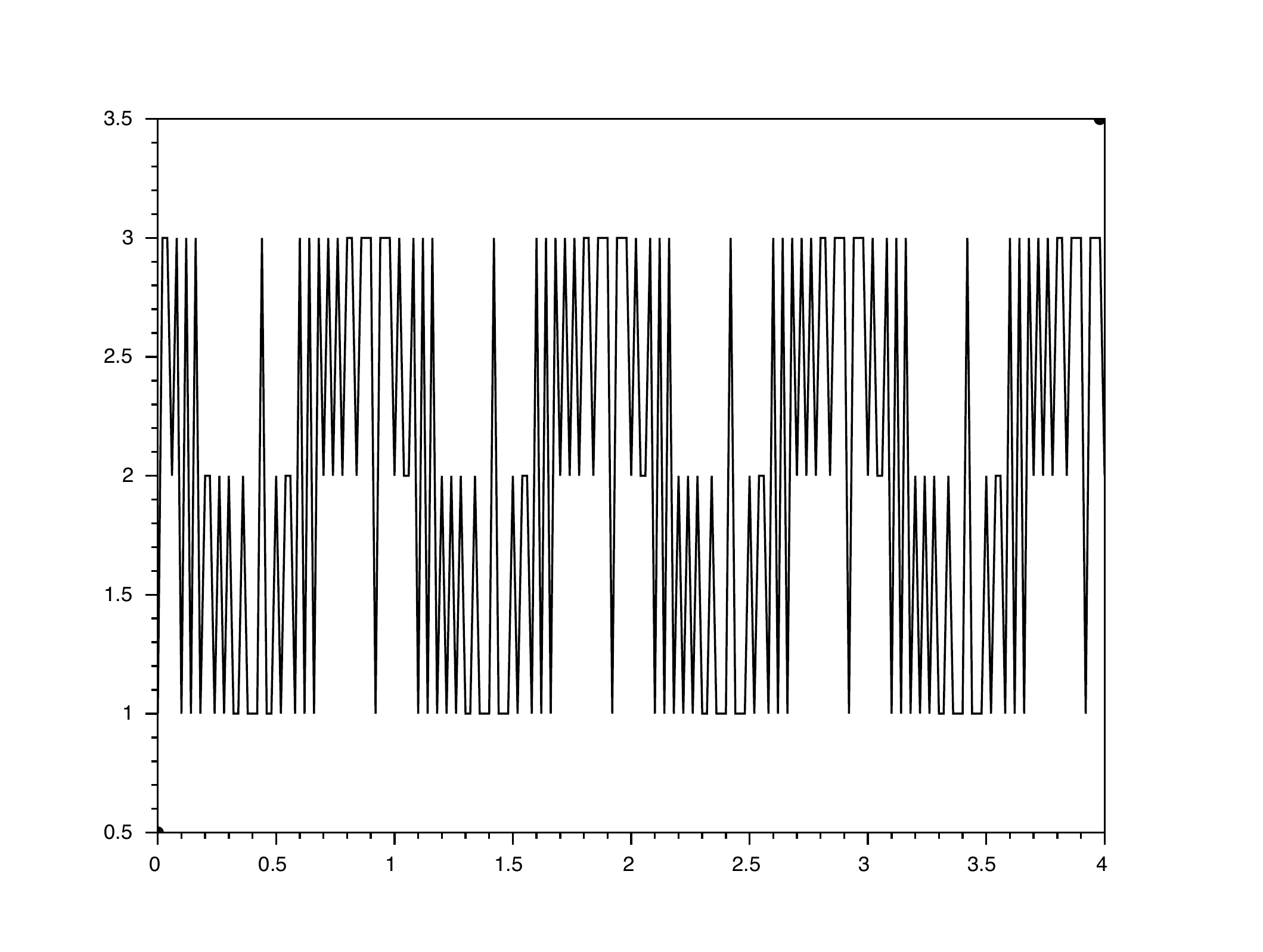}
    \caption{Value function and optimal switching for the DC/AC inverter.}
    \label{fig:inverter1}
    \includegraphics[width=.49\textwidth]{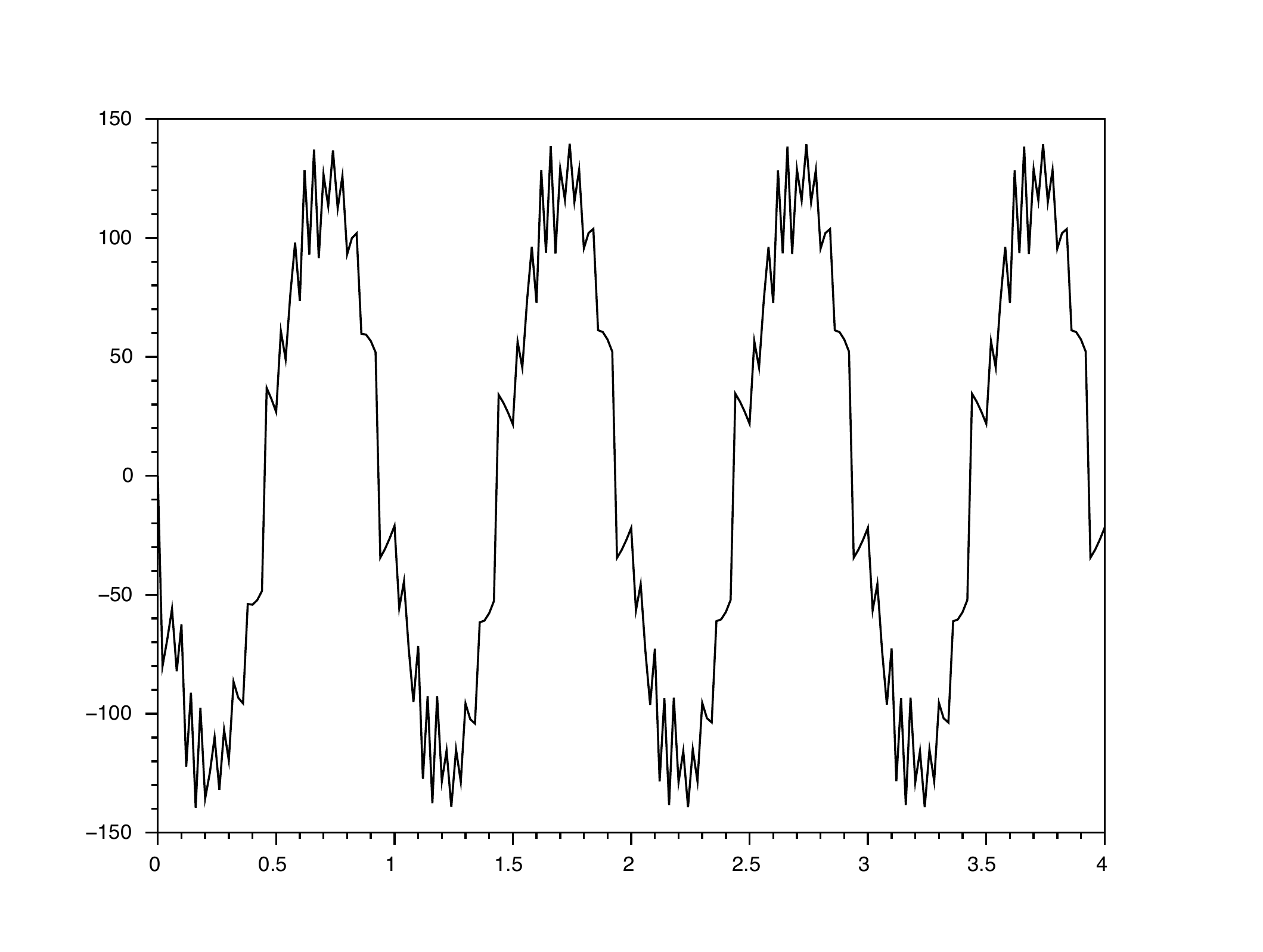}
    \includegraphics[width=.49\textwidth]{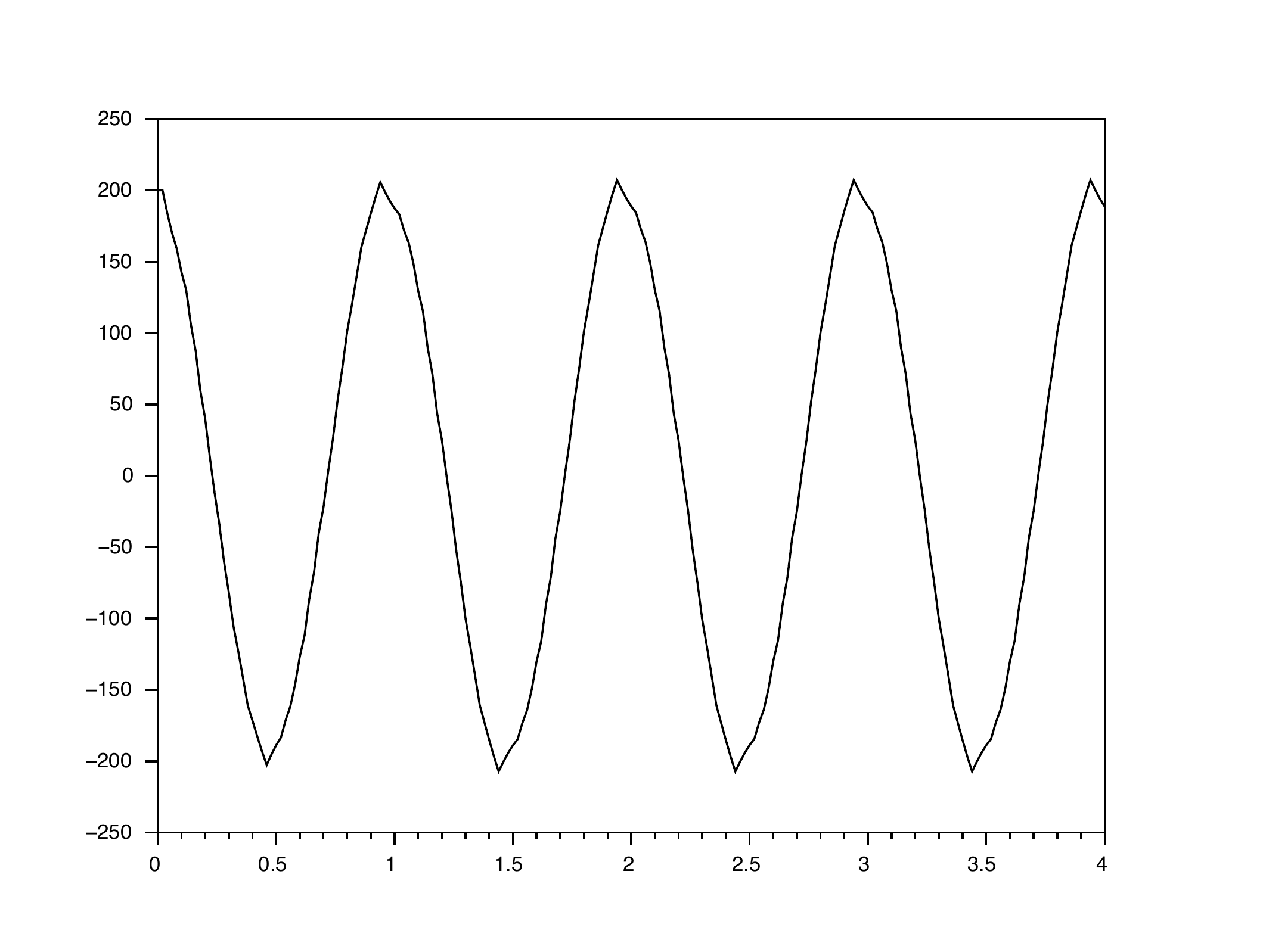}
    \caption{Trajectories $X_1(t)$ and $X_2(t)$ for the DC/AC inverter.}
    \label{fig:inverter2}
    \includegraphics[width=.7\textwidth]{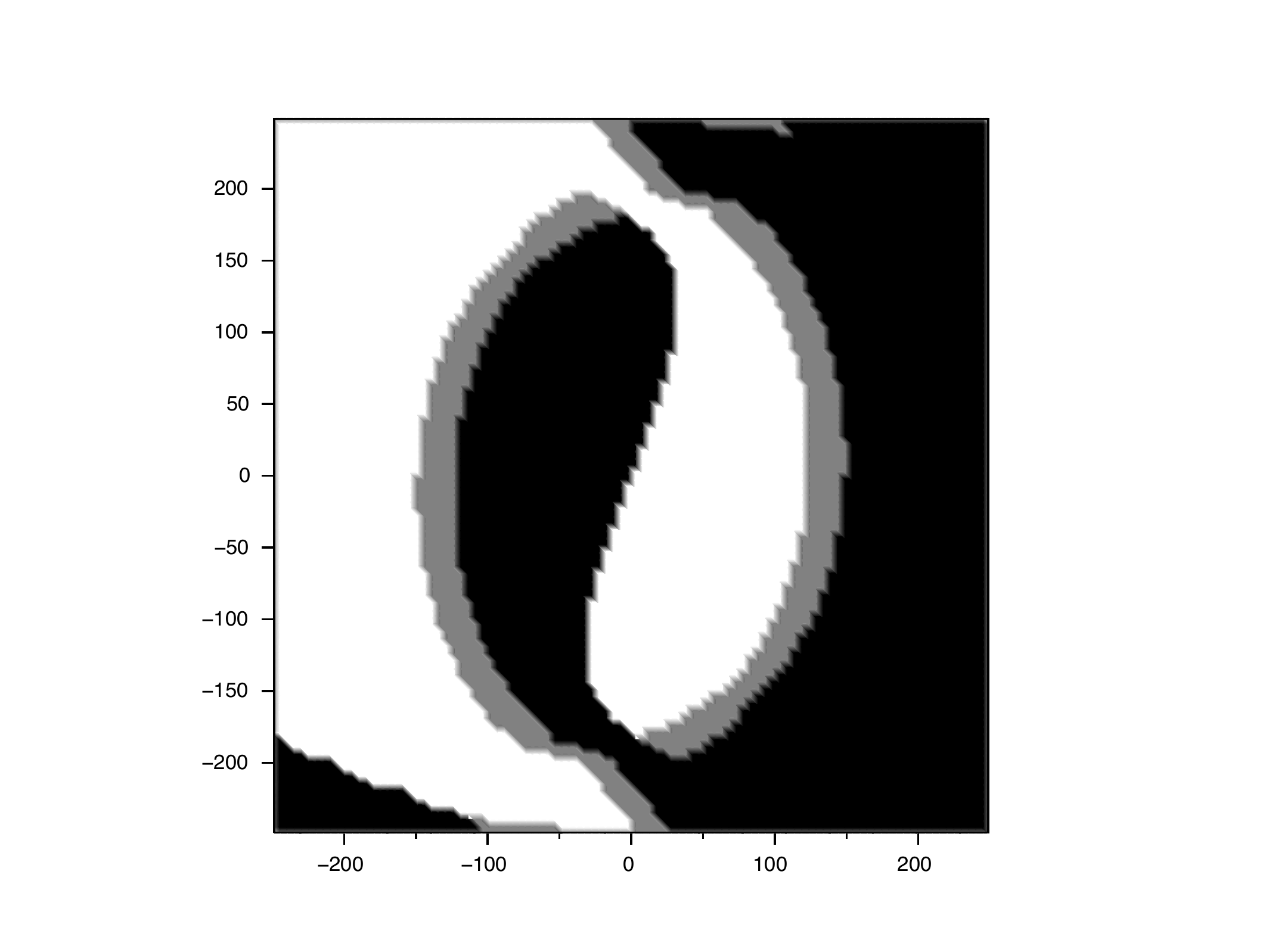}
    \caption{Optimal switching map for $q=2$ for the DC/AC inverter.}
    \label{fig:inverter3}
\end{figure}

In this device, a DC source generates an AC output by means of a suitable operation of the switches $S_1,\ldots,S_4$, as well as a suitable choice of the three components ($R$, $L$ and $C$) which appear in series in the $RLC$ load. Following \cite{ChSan}, we consider as state variables $X_1=i_L$ (the current through the inductor $L$, i.e., through the load) and $X_2=v_C$ (the voltage across the capacitor $C$), the state equations being
\begin{equation}\label{inverter_model}
\begin{cases}
\displaystyle\dot X_1(t) = \frac{V_{DC}}{L}(Q(t)-2) - \frac{R}{L} X_1(t) - \frac{1}{L} X_2(t) \\
\displaystyle\dot X_2(t) = \frac{1}{C} X_1(t).
\end{cases}
\end{equation}
The physical meaning of the discrete state variable depends on the state of the switches $S_1,\ldots,S_4$, and more precisely
\[
Q(t) =
\begin{cases}
1 & \text{if $S_1,S_3=OFF$ and $S_2,S_4=ON$} \\
2 & \text{if $S_1,S_4=OFF$ and $S_2,S_3=ON$} \\
3 & \text{if $S_2,S_4=OFF$ and $S_1,S_3=ON$.} 
\end{cases}
\]
The cost functional is defined so as to force the system to evolve (approximately) along an ellipse of the state space (see \cite{ChSan}), namely
\[
\frac{x_1^2}{a^2} + \frac{x_2^2}{b^2} = c,
\]
in which the constants $a$ and $b$ are defined in terms of the physical parameters $R$, $L$, $C$ and of the desired pulsation $\omega$. This makes it natural to define the running cost as
\begin{equation}\label{cost_inverter}
\ell(x,q,\alpha) = \left(\frac{x_1^2}{a^2} + \frac{x_2^2}{b^2} - c\right)^2.
\end{equation}
With the parameters chosen (see Table \ref{tab:parameter3}), $a\approx 0.84$, $b\approx 1.34$ and the required output of the system would be given by two sinusoids of amplitude respectively 126 A for $X_1$ and 200 V for $X_2$, both at the frequency of 1 Hz. The approximate solution has been computed on a $100\times 100$ grid on the domain $[-250,250]^2$, with $\Delta t=0.01$, and state constraint boundary conditions have been treated by penalization, assigning a stopping cost of $5\cdot 10^8$ on the boundary. The effect of the lack of full controllability is apparent in Fig. \ref{fig:inverter1}, which shows one component of the value function (they are practically undistinguishable from one another) and the optimal switching with respect to time. Fig. \ref{fig:inverter2} shows the output $(X_1(t),X_2(t))$ of the controlled system, whereas, as an example, Fig. \ref{fig:inverter3} reports the switching map of the second component of the state space. Here, the optimal solution is to keep $Q(t)=2$ in grey regions, commute to $Q(t)=1$ in black regions and to $Q(t)=3$ in white regions.

Finally, Table \ref{tab:inverter} compares the two numerical solvers (VI and MPI). In this last test, the stopping condition has been computed on the relative $l^1$ update,
$$
\frac{\|\bm v_j - \bm v_{j-1}\|_1}{\|\bm v_j\|_1} < \epsilon,
$$
to avoid problems with both high values of the solution and the occurrence of a discontinuity caused by the lack of controllability.

\begin{table}
\centering
\begin{tabular}{|c|c|c|c|c|c|c|c|c|}
\hline
$V_{DC}$ & $R$ & $L$ & $C$ & $\omega$ & $c$ & $\lambda$ & $x_1$ & $x_2$ \\
\hline \hline
$200$ [V] & $0.7$ [$\Omega$] & $0.1$ [H] & $0.1$ [F] & $2\pi$ [s$^{-1}$] & $22500$ & $1$ & $0$ & $200$ \\
\hline
\end{tabular}
\caption{Choice of parameters, inverter}\label{tab:parameter3}
\bigskip
\begin{tabular}{|c|c|c|c|}
\hline
$\epsilon$ & $N_V$ & $N_P$ \\
\hline \hline
$10^{-3}$ & $469$ & $469$  \\
\hline
$10^{-6}$ & $13481$ & $13491$  \\
\hline
\end{tabular}
\caption{Number of iterations for a given tolerance $\epsilon$, DC/AC inverter test.}
\label{tab:inverter}
\end{table}

\section*{Conclusions}

We have constructed and validated a Semi-Lagrangian scheme for hybrid Dynamic Programming problems in infinite horizon form. The numerical scheme has been made more efficient by a Policy Iteration type solver. Numerical tests performed on examples of varying complexity show that the scheme is robust and that the approximate optimal control policy obtained is stable and accurate, although the complexity remains critical with respect to the dimension of the state space.

This validation suggests that this could be a feasible method to design optimization-based static controllers in low dimension.

\end{document}